\documentclass{amsart}%
\usepackage{amsfonts}
\usepackage{amsmath}
\usepackage{amssymb}
\usepackage{amsthm}
\usepackage{booktabs}
\usepackage{mathtools}
\usepackage[dvipsnames]{xcolor}
\usepackage[colorlinks, linkcolor=BrickRed,citecolor=Violet]{hyperref}
\usepackage[noabbrev,nameinlink]{cleveref}
\usepackage[shortlabels]{enumitem}
\usepackage{tikz}
\usepackage{tikz-cd}
\usetikzlibrary{calc, shapes,arrows,positioning}
\tikzset{>=stealth',
  head/.style = {fill = white, text=black},
  plaque/.style = {draw, rectangle, minimum size = 10mm}, 
  pil/.style={->,thick},
  junct/.style = {draw,circle,inner sep=0.5pt,outer sep=0pt, fill=black}
  }
\usepackage{varwidth}
\usepackage{subcaption}
\usepackage{multirow}
\usepackage{rotating}
\usepackage[boxsize=1em]{ytableau}

\makeatletter
\newtheorem*{rep@theorem}{\rep@title}
\newcommand{\newreptheorem}[2]{%
\newenvironment{rep#1}[1]{%
 \def\rep@title{#2 \ref{##1}}%
 \begin{rep@theorem}}%
 {\end{rep@theorem}}}
\makeatother

\makeatletter
\newtheorem*{rep@corollary}{\rep@title}
\newcommand{\newrepcorollary}[2]{%
\newenvironment{rep#1}[1]{%
 \def\rep@title{#2 \ref{##1}}%
 \begin{rep@corollary}}%
 {\end{rep@corollary}}}
\makeatother

\newtheorem{theorem}{Theorem}[section]
\newreptheorem{theorem}{Theorem}

\newtheorem{lemma}[theorem]{Lemma}
\newreptheorem{corollary}{Corollary}

\newtheorem{proposition}[theorem]{Proposition}

\theoremstyle{definition}

\numberwithin{equation}{section}

\newenvironment{cd}[1][]
  {\begin{center}\begin{tikzcd}[#1]}
  {\end{tikzcd}\end{center}}

\newcommand{\bZ}{\mathbb{Z}}
\newcommand{\cE}{\mathcal{E}}

\DeclareMathOperator{\Inc}{Inc}
\DeclareMathOperator{\KPro}{K-Pro}
\DeclareMathOperator{\maj}{maj}
\DeclareMathOperator{\Packed}{Pack}
\newcommand{\IncP}{\Inc_{\Packed}}

\DeclareMathOperator{\Row}{Row}
\DeclareMathOperator{\SYT}{SYT}

\newcommand{\qbinom}{\genfrac{[}{]}{0pt}{}}

\newcommand{\blank}{\phantom{2}}

\title[Curious cyclic sieving on increasing tableaux]{Curious cyclic sieving \\ on increasing tableaux}
\author[Gaetz]{Christian Gaetz}
\address{Department of Mathematics, Harvard University, Cambridge, MA, 02138, USA}
\email{gaetz@math.harvard.edu}

\author[Pechenik]{Oliver Pechenik}
\address{Department of Combinatorics \& Optimization, University of Waterloo, Waterloo, ON, N2L 3G1, Canada}
\email{oliver.pechenik@uwaterloo.ca}

\author[Striker]{Jessica Striker}
\address{Department of Mathematics, North Dakota State University, Fargo, ND, 58102, USA}
\email{jessica.striker@ndsu.edu}

\author[Swanson]{Joshua P. Swanson}
\address{Department of Mathematics, University of Southern California, Los Angeles, CA, 90007, USA}
\email{swansonj@usc.edu}

\date{\today}

\begin{document}

\begin{abstract}
  We prove a cyclic sieving result for the set of $3 \times k$ packed increasing tableaux with maximum entry $m \coloneqq 3+k$ under K-promotion. The ``curiosity'' is that the sieving polynomial arises from the $q$-hook formula for standard tableaux of ``toothbrush shape'' $(2^3, 1^{k-2})$ with $m+1$ boxes, whereas K-promotion here only has order $m$.
\end{abstract}
\maketitle

\section{Introduction}

\subsection{Increasing tableaux and main result}

An \textit{increasing tableau} is a filling $T$ of the diagram of an \textit{integer partition} $\lambda$ with positive integers that strictly increase along rows and columns; $T=\ytableaushort{125,23}$  is an increasing tableau of shape $\lambda=(3,2)$. Let $\Inc^m(\lambda)$ be the set of increasing tableaux with maximum entry $\max(T)$ at most $m$. We call $T$ \textit{packed} if each value $1, 2, \ldots, \max(T)$ appears at least once. Let $\IncP^m(\lambda)$ be the set of packed increasing tableaux of shape $\lambda$ and maximum entry exactly $m$; \ytableaushort{124,23} \ is an increasing tableau in $\IncP^4(3,2)$. Let $\KPro$ denote the \textit{K-promotion} operator on $\Inc^m(\lambda)$; see \Cref{sec:back:k_pro} for details. K-promotion was introduced by the second author in \cite{MR3207480} building on work of Thomas--Yong \cite{MR2491941}, and it has been further studied in \cite{MR3537922,MR3603321,MR3711092,MR4190062}. Finally, let $\SYT(\lambda) \coloneqq \IncP^{|\lambda|}(\lambda)$ denote the set of all \textit{standard Young tableaux} of shape $\lambda$ (see e.g.~\cite[\S7.10]{MR1676282}).

\textit{Dynamical algebraic combinatorics} is concerned with the properties of explicit combinatorial discrete dynamical systems, such as the number and sizes of orbits for a bijection applied iteratively to a finite set \cite{MR3526426,MR3585535}. The promotion operator $\KPro$ on $\Inc^{a+b}(a \times b)$ has order $a+b$ \cite[Cor.~4.10]{MR3603321}, where $a \times b$ denotes the rectangular partition $(b^a)$ with $a$ rows of length $b$. The \textit{cyclic sieving phenomenon} (\textit{CSP}) of Reiner--Stanton--White \cite{MR2087303} encodes the orbit structure of a cyclic action on a finite set in evaluations of a polynomial at roots of unity; see \Cref{sec:back:CSP}.

We consider the very special case of $\KPro$ acting on $\IncP^{3+k}(3 \times k)$, certain rectangular increasing tableaux with three rows. Correspondingly, it turns out we will be interested in the ``toothbrush-shaped'' standard Young tableaux $\SYT(2^3, 1^{k-2})$; see \Cref{fig:toothbrush}. Perhaps surprisingly, these sets are equinumerous.

\begin{proposition}\label{prop:main}
For all $k>1$, we have
   \[ \left|\IncP^{3+k}(3 \times k)\right| = \left|\SYT(2^3, 1^{k-2})\right|. \]
\end{proposition}

While we are able to give explicit bijections between the sets of \Cref{prop:main}, we have not been able to identify any canonical structure-preserving bijection. See \Cref{sec:bijections} for further discussion.

\begin{figure}[htb]
\ytableaushort{14,26,3{10},5,7,8,9}
\caption{A standard Young tableau of toothbrush shape $(2^3,1^{k-2})$ with $k=6$.}\label{fig:toothbrush}
\end{figure}

For $\lambda \vdash N$, let
\[
f^\lambda(q) \coloneqq  \frac{[N]_q!}{\prod_{c \in \lambda} [h_c]_q}
\]
be a $q$-analogue of the hook length formula; see \Cref{sec:back:hook_length}. Our main result is as follows.

\begin{theorem}\label{thm:main}
  The triple 
  \[
  \left(\IncP^{3+k}(3 \times k), \langle \KPro\rangle, f^{(2^3, 1^{k-2})}(q)\right)
  \]
  exhibits the cyclic sieving phenomenon, and $\KPro^{3+k} = 1$.
\end{theorem}

We note that the cyclic sieving polynomial in \Cref{thm:main} is the $q$-hook length formula for the toothbrush shape $(2^3, 1^{k-2})$, while the tableaux appearing are instead increasing tableaux of rectangular shape $3 \times k$.

\subsection{Comparison to existing results}\label{sec:existing}

Several results similar to \Cref{prop:main} and \Cref{thm:main} have appeared for increasing tableaux with special shapes and contents. There is as yet no unifying generalization, and finding further examples would be of interest.

\begin{itemize}
  \item[(CSP.1)]\label{thm:rhoades} Rhoades \cite[Thm.~1.3]{MR2557880} showed that
    \[ (\IncP^{ab}(a \times b), \langle\KPro\rangle, f^{a \times b}(q)) \]
  exhibits the CSP, where $\KPro$ has order $ab$. Note that in this case, the tableaux are standard Young tableaux and the cyclic sieving polynomial is the $q$-hook length formula for the shape of the tableaux.

  \item[(CSP.2)]\label{thm:pechenik} The second author \cite[Thm.~1.2]{MR3207480} showed that
    \[ (\IncP^m(2 \times k), \langle \KPro\rangle, f^{(m-k, m-k, 1^{2k-m})}(q)) \]
  exhibits the CSP, where $\KPro$ has order $m$. Compared with (CSP.1), the tableaux here are of more restrictive shape but more general maximum entry. Note that, as in \Cref{thm:main}, the cyclic sieving polynomial here is the $q$-hook length formula for the shape $(m-k, m-k, 1^{2k-m})$ which different from the shape of the tableaux in question.
  
  The second author, moreover, gave an explicit bijection \cite[Thm.~1.1]{MR3207480} between two-row rectangular increasing tableaux $\IncP^m(2 \times k)$ and standard Young tableaux of ``pennant shape'' $(m-k, m-k, 1^{2k-m})$. This bijection is not $\KPro$-equivariant, but is equivariant for the related involution \emph{K-evacuation} and also preserves the \emph{descents} of the tableaux.
  
  When $m=2+k$, we have direct two-row analogues of \Cref{prop:main} and \Cref{thm:main} arising from
    \[ \left|\IncP^{2+k}(2 \times k)\right| = \left|\SYT(2^2, 1^{k-2})\right|. \]

  \item[(CSP.3)]\label{thm:psv} Pressey--Stokke--Visentin \cite[Thm.~3.7]{MR3537922} showed that
    \[ (\IncP^{m}(r, 1^s), \langle \KPro\rangle, f^{(m-s,1^s)}(q) f^{(m-r+1, 1^{r+s-m})}(q)) \]
  exhibits the CSP, where $\KPro$ has order $m-1$. We note that in this case the order $m-1$ of K-promotion differs from the number $m+s+1$ of cells in the pairs of standard tableaux, as well as from the maximum entry $m$ of the increasing tableaux.
\end{itemize}

\subsection{Potential generalizations}
In light of the results (CSP.1), (CSP.2), and (CSP.3), which all give cyclic sieving phenomena for K-promotion on various sets of increasing tableaux with sieving polynomial a product of $q$-hook polynomials, one might ask if Theorem~\ref{thm:main} could be generalized while maintaining this property. However the most natural extensions of Theorem~\ref{thm:main} in this direction, to $\IncP^{4 + k}(4 \times k)$ or to $\IncP^{(3+k)+1}(3 \times k)$, will not work. One can compute that:
\begin{align*}
\left| \IncP^{4+4}(4 \times 4) \right| &= 2 \cdot 31 \\
\left| \IncP^{(3+7)+1}(3 \times 7) \right| &= 5 \cdot 11 \cdot 67.
\end{align*}
The value $f^{\lambda}(1)$ has largest prime divisor at most $N$, the number of boxes of $\lambda$, so the prime factors of 31 and 67 above severely restrict which $q$-hook polynomials could appear as factors in a potential sieving polynomial, and it can easily be verified that none of the possibilities is in fact a sieving polynomial for K-promotion on these sets of increasing tableaux. Note also that, for general $a,b,m \in \mathbb{N}$, the order of K-promotion on $\IncP^m(a \times b)$ is strictly greater than $m$ and moreover the order is unknown (cf.\ \cite{MR3207480,MR4190062}).

\subsection{Organization}
The rest of the paper is organized as follows. In \Cref{sec:back}, we give background on K-promotion, the cyclic sieving phenomenon, hook length formulas, and rowmotion on order ideals. In \Cref{sec:proof}, we prove \Cref{prop:main} and \Cref{thm:main}. In \Cref{sec:bijections}, we discuss bijectivity.

\section{Background}\label{sec:back}

\subsection{K-promotion}\label{sec:back:k_pro}

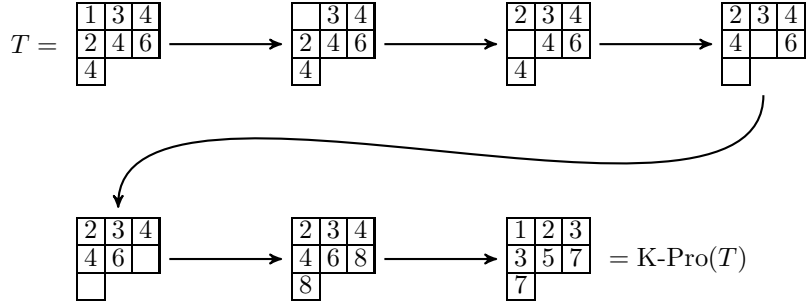
\begin{figure}[htb]
  \centering

\begin{tikzpicture}
\node (A) {\ytableaushort{134,246,4}};
\node[left = 0 of A] (T) {$T=$};
\node[right=1.5 of A] (B) {\ytableaushort{\blank 34, 246,4}};
\node[right=1.5 of B] (C) {\ytableaushort{2 3 4, \blank 46,4}};
\node[right=1.5 of C] (D) {\ytableaushort{234, 4 \blank 6, \blank}};
\node[below=1.5 of A] (E) {\ytableaushort{234, 46 \blank, \blank}};
\node[right=1.5 of E] (F) {\ytableaushort{234,468,8}};
\node[right=1.5 of F] (G) {\ytableaushort{123, 357,7}};
\node[right=0.0 of G] (H) {$=\KPro(T)$};
\path (A) edge[pil]  (B);
\path (B) edge[pil]  (C);
\path (C) edge[pil]  (D);
\draw[->, thick] (D.south) .. controls ([yshift=-3cm] D) and ([yshift=3cm] E) .. (E.north);
\path (E) edge[pil]  (F);
\path (F) edge[pil]  (G);
\end{tikzpicture}

  \caption{An example of K-promotion for $T \in \Inc^7(3,3,1)$.}
  \label{fig:KPro}
\end{figure}

Thomas--Yong \cite{MR2491941} introduced K-jeu de taquin for increasing tableaux. K-promotion on $\Inc^m(\lambda)$ was built out of sliding moves in \cite{MR3207480} as follows. The \textit{southeast neighbors} of a cell are the (at most two) adjacent cells immediately south or east of it; see \Cref{fig:KPro}. Let $T \in \Inc^m(\lambda)$ be an increasing tableau\footnote{In contrast to \cite{MR3207480}, we do not require $T$ to be packed, and the value $m$ is not required to appear in~$T$.}. Delete the entry $1$ from $T$, leaving an empty cell. Repeatedly perform the following operation simultaneously on all empty cells until no empty cell has a southeast neighbor. Label each empty cell by the minimal label of its southeast neighbor(s) and then remove that label from the southeast neighbor(s) in which it appears. If an empty cell has no southeast neighbors, it remains unchanged. Finally we obtain $\KPro(T)$ by labeling all empty cells by $m+1$ and then subtracting $1$ from every label.

\subsection{The cyclic sieving phenomenon}\label{sec:back:CSP}

Let $X$ be a finite set on which a cyclic group $C$ acts, and let $f(q) \in \bZ_{\geq 0}[q]$ be a polynomial. We say the triple $(X, C, f(q))$ exhibits the \textit{cyclic sieving phenomenon} (\textit{CSP}) if the number of fixed points of any element $\sigma \in C$ of order $d$ is $f(\exp(2\pi i/d))$ \cite{MR2087303}. In particular, $|X| = f(1)$, so \Cref{prop:main} follows from \Cref{thm:main}. See \cite{MR2866734} for a nice survey article on the CSP.

\subsection{Hook lengths and \texorpdfstring{$q$}{q}-analogues}\label{sec:back:hook_length}

The \textit{hook length formula} is 
\[
|\SYT(\lambda)| = \frac{N!}{\prod_{c \in \lambda} h_c},
\]
where $\lambda$ has $N$ boxes and $h_c$ is the \textit{hook length} of the cell $c$ \cite[p.373]{MR1676282}. There is a natural $q$-analogue of the hook length formula which enumerates standard Young tableaux by their \textit{major index}:
  \[ \sum_{T \in \SYT(\lambda)} q^{\maj(T)} = q^{b(\lambda)} \frac{[N]_q!}{\prod_{c \in \lambda} [h_c]_q} = q^{b(\lambda)} f^\lambda(q), \]
  where $[d]_q \coloneqq \frac{1-q^d}{1-q}$, $[N]_q! \coloneqq [N]_q [N-1]_q \cdots [1]_q$, and $b(\lambda) \coloneqq \sum_{i \geq 1} (i-1) \lambda_i$ \cite[Cor.~7.21.5]{MR1676282}. The special case $\lambda = (a+1, 1^b)$ yields the \textit{$q$-binomial coefficients}
    \[ f^{(a+1, 1^b)}(q) = \frac{[a+b]_q!}{[a]_q! [b]_q!} \eqqcolon \qbinom{a+b}{a}_q. \]

\subsection{Order ideals and rowmotion}
Let $[a] \times [b]$ denote the poset which is the product of two chains, i.e.\ the grid $\{(i, j) : 1 \leq i \leq a, 1 \leq j \leq b\}$ under the component-wise partial order. Let $J(P)$ denote the set of \textit{order ideals} of a poset $P$, which for $[a] \times [b]$ will be identified with lower-right-justified Young diagrams in the $a \times b$ rectangle; see \Cref{fig:lattice_path}. 

\textit{Rowmotion} on $P$ is the bijection $\Row \colon J(P) \to J(P)$ defined by sending $J \in J(P)$ to the order ideal generated by the minimal elements of the complement $P - J$ \cite{MR0349497,MR2950491}. There is an alternative characterization in terms of \emph{toggles} \cite{MR1356845}, which was used in \cite{MR3603321} to give an equivariant bijection between K-promotion on increasing tableaux (using an alternative interpretation in terms of \emph{K-Bender-Knuth involutions}) and rowmotion on order ideals of related posets. We give a self-contained characterization of this bijection specialized to our context in the next section.

\section{Proof of the main results}\label{sec:proof}
    We begin by describing an equivariant bijection between $\Inc^{a+b}(a \times b)$ and order ideals $J([a] \times [b])$. Suppose $T \in \Inc^{a+b}(a \times b)$. Index the cells of $T$ with $(1, 1)$ at the upper left and $(a, b)$ at the lower right. Let $J = \{(i, j) : T(i, j) = i+j\}$. Consider traveling from $(i, j) \in J$ to $(a, b)$ using adjacent horizontal or vertical steps. At each step the value of $i+j$ increases by exactly one and the value of $T$ increases by at least one, so
      \[ a+b \geq T(a,b) \geq T(i, j) + (a-i) + (b-j) = a+b. \]
    Thus $T$ increases by exactly one at every step, so everything weakly southeast of $(i, j)$ belongs to $J$, and $J$ is an order filter of $[a]\times[b]$. By the same argument, $K = \{(i, j) : T(i, j) = i+j-1\}$ is $([a] \times [b]) - J$ and is an order ideal. It is easy to see that $T \mapsto (J, K)$ is bijective. Hence, the map $\Theta : \Inc^{a+b}(a \times b) \to J([a] \times [b])$ defined by $T \mapsto \{(a-i+1, b-i+1) : (i, j) \in J\}$ is a bijection.

\begin{lemma}\label{lem:lattice_path}
  For any positive integers $a$ and $b$, the bijection $\Theta$ intertwines rowmotion with K-promotion:
  \begin{cd}
    \Inc^{a+b}(a \times b) \rar{\Theta} \dar[swap]{\KPro}
      & J([a] \times [b]) \dar{\Row} \\
    \Inc^{a+b}(a \times b) \rar{\Theta}
      & J([a] \times [b])
  \end{cd}
  \begin{proof}
   We have already observed that $\Theta$ is bijective. For equivariance, consider the effect of $\KPro$ on $(J, K)$ and suppose that $\KPro(T) \mapsto (J', K')$. Each rectangle weakly northwest of a maximal element $(i, j)$ of $K$ is completely determined, and K-promotion fixes that rectangle except for $(i, j)$. Since the east and/or south neighbors of $(i, j)$ belong to $J$ and so are larger by $2$, K-promotion increases $T(i, j)$ by $1$. This holds even when $(i, j) = (a, b)$, so $(i, j) \in J'$. Hence elements weakly southeast of maximal elements of $K$ belong to $J'$ and elements weakly northwest, not including the element itself, belong to $K'$. All that remains are elements strictly northeast or southwest of maximal elements of $K$, necessarily in $J$. It is not difficult to see that these are decremented by $1$ by K-promotion and hence belong to $K'$. Thus $J'$ is precisely the order filter generated by the maximal elements of the complement of $J$, and equivariance follows.
  \end{proof}
\end{lemma}

The existence of a bijection making the diagram of \Cref{lem:lattice_path} commute was previously established in \cite{MR3603321} in more generality, but without making the map explicit.

\begin{figure}[ht]
  \centering

\begin{tikzpicture}[inner sep=0in,outer sep=0in]
\node (n) {\begin{varwidth}{5cm}{
\ytableausetup{boxsize=1.45em}
\begin{ytableau}
1 & 2 & *(gray!40)4 & *(gray!40)5 & *(gray!40)6 \\
2 & 3 & *(gray!40)5 & *(gray!40)6 & *(gray!40)7 \\
3 & *(gray!40)5 & *(gray!40)6 & *(gray!40)7 & *(gray!40)8 \\
*(gray!40)5 & *(gray!40)6 & *(gray!40)7 & *(gray!40)8 & *(gray!40)9
\end{ytableau}}\end{varwidth}};
\ytableausetup{boxsize=1em}
\draw[very thick,blue] (n.south west)--++(0,1.5em)--++(1.5em,0)--++(0,1.5em)--++(1.5em,0)--++(0,1.5em)--++(0,1.5em)--++(1.5em,0)--++(1.5em,0)--++(1.5em,0);
\end{tikzpicture}
  \caption{A rectangular increasing tableau with its corresponding Young diagram and upper order ideal $J$ from \Cref{lem:lattice_path} shaded. The edge of the Young diagram is determined by the horizontally or vertically adjacent cells whose entries differ by exactly $2$.}
  \label{fig:lattice_path}
\end{figure}

\begin{lemma}\label{cor:Inc_CSP}
  For any positive integers $a$ and $b$, the triple 
  \[
  \left(\Inc^{a+b}(a \times b), \langle \KPro \rangle, \qbinom{a+b}{a}_q\right)
  \]
  exhibits the cyclic sieving phenomenon.

  \begin{proof}
    Stanley gave an equivariant bijection between $J([a] \times [b])$ under rowmotion and skew standard Young tableaux with two rows of lengths $a$ and $b$ that do not overlap (denoted $\SYT((a+b, b)/(b))$) under promotion \cite[p.8]{MR2515772}. In \cite[\S3.1]{MR2950491}, the third author and Williams coined the term ``rowmotion,'' noted that $\SYT((a+b, b)/(b))$ under promotion is in equivariant bijection with $\binom{[a+b]}{b}$ under cyclic rotation, and noted that cyclic sieving applies by Reiner--Stanton--White's foundational example \cite[Thm.~1.1]{MR2087303}. The result follows by combining these observations with \Cref{lem:lattice_path}.
  \end{proof}
    \end{lemma}

In general, the relation between the K-promotion orbits of $\Inc^m(\lambda)$ and those of the subset $\IncP^m(\lambda)$ is fairly complicated (cf.\ \cite[Theorem~6.1]{MR3846229}). However, in the following case of particular relevance to this paper, it is much simpler.

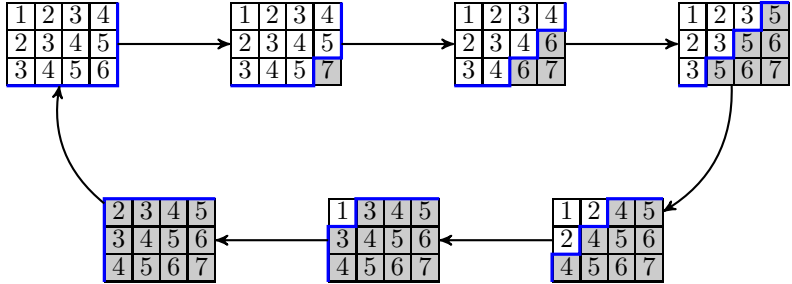
\begin{figure}[ht]
  \centering

\begin{tikzpicture}[inner sep=0in,outer sep=0in]
\node (A) {
\begin{varwidth}{5cm}{
\begin{ytableau}
1 & 2 & 3 & 4 \\
2 & 3 & 4 & 5 \\
3 & 4 & 5 & 6
\end{ytableau}}\end{varwidth}
};
\draw[very thick,blue] (A.south west)--++(1.05em,0)--++(1.05em,0)--++(1.05em,0)--++(1.05em,0)--++(0,1.05em)--++(0,1.05em)--++(0,1.05em);

\node[right=1.5 of A] (B) {
\begin{varwidth}{5cm}{
\begin{ytableau}
1 & 2 & 3 & 4 \\
2 & 3 & 4 & 5 \\
3 & 4 & 5 & *(gray!40)7
\end{ytableau}}\end{varwidth}
};
\draw[very thick,blue] (B.south west)--++(1.05em,0)--++(1.05em,0)--++(1.05em,0)--++(0,1.05em)--++(1.05em,0)--++(0,1.05em)--++(0,1.05em);

\node[right=1.5 of B] (C) {
\begin{varwidth}{5cm}{
\begin{ytableau}
1 & 2 & 3 & 4 \\
2 & 3 & 4 & *(gray!40)6 \\
3 & 4 & *(gray!40)6 & *(gray!40)7
\end{ytableau}}\end{varwidth}
};
\draw[very thick,blue] (C.south west)--++(1.05em,0)--++(1.05em,0)--++(0,1.05em)--++(1.05em,0)--++(0,1.05em)--++(1.05em,0)--++(0,1.05em);

\node[right=1.5 of C] (D) {
\begin{varwidth}{5cm}{
\begin{ytableau}
1 & 2 & 3 & *(gray!40)5 \\
2 & 3 & *(gray!40)5 & *(gray!40)6 \\
3 & *(gray!40)5 & *(gray!40)6 & *(gray!40)7
\end{ytableau}}\end{varwidth}
};
\draw[very thick,blue] (D.south west)--++(1.05em,0)--++(0,1.05em)--++(1.05em,0)--++(0,1.05em)--++(1.05em,0)--++(0,1.05em)--++(1.05em,0);

\node[below left=1.5 and 0.2 of D] (E) {
\begin{varwidth}{5cm}{
\begin{ytableau}
1 & 2 & *(gray!40)4 & *(gray!40)5 \\
2 & *(gray!40)4 & *(gray!40)5 & *(gray!40)6 \\
*(gray!40)4 & *(gray!40)5 & *(gray!40)6 & *(gray!40)7
\end{ytableau}}\end{varwidth}
};
\draw[very thick,blue] (E.south west)--++(0,1.05em)--++(1.05em,0)--++(0,1.05em)--++(1.05em,0)--++(0,1.05em)--++(1.05em,0)--++(1.05em,0);

\node[left=1.5 of E] (F) {
\begin{varwidth}{5cm}{
\begin{ytableau}
1 & *(gray!40)3 & *(gray!40)4 & *(gray!40)5 \\
*(gray!40)3 & *(gray!40)4 & *(gray!40)5 & *(gray!40)6 \\
*(gray!40)4 & *(gray!40)5 & *(gray!40)6 & *(gray!40)7
\end{ytableau}}\end{varwidth}
};
\draw[very thick,blue] (F.south west)--++(0,1.05em)--++(0,1.05em)--++(1.05em,0)--++(0,1.05em)--++(1.05em,0)--++(1.05em,0)--++(1.05em,0);

\node[left=1.5 of F] (G) {
\begin{varwidth}{5cm}{
\begin{ytableau}
*(gray!40)2 & *(gray!40)3 & *(gray!40)4 & *(gray!40)5 \\
*(gray!40)3 & *(gray!40)4 & *(gray!40)5 & *(gray!40)6 \\
*(gray!40)4 & *(gray!40)5 & *(gray!40)6 & *(gray!40)7
\end{ytableau}}\end{varwidth}
};
\draw[very thick,blue] (G.south west)--++(0,1.05em)--++(0,1.05em)--++(0,1.05em)--++(1.05em,0)--++(1.05em,0)--++(1.05em,0)--++(1.05em,0);

\path (A) edge[pil]  (B);
\path (B) edge[pil]  (C);
\path (C) edge[pil]  (D);
\path (D) edge[pil,bend left=30]  (E);
\path (E) edge[pil]  (F);
\path (F) edge[pil]  (G);
\path (G) edge[pil,bend left=30]  (A);
\end{tikzpicture}

  \caption{This seven-cycle is the exceptional orbit of K-promotion $\cE = \Inc^7(3 \times 4) - \IncP^7(3 \times 4),$ as described by \Cref{lem:decomp}.}
  \label{fig:exceptional_orbit}
\end{figure}

\begin{lemma}\label{lem:decomp}
  We have
    \[ \Inc^{a+b}(a \times b) = \IncP^{a+b}(a \times b) \sqcup \cE, \]
  where $\cE$ is the orbit of the tableau $T(i, j) = i+j-1$ under $\KPro$, and $|\cE| = a+b$.

  \begin{proof}
    This is easy to see from \Cref{lem:lattice_path} and \Cref{fig:exceptional_orbit}. In the notation of the proof of \Cref{lem:lattice_path}, the packed tableaux in $\Inc^{a+b}(a \times b)$ are precisely those for which some anti-diagonal $\{(i, j) : i+j=c\}$ has elements of both $J$ and $K$.
  \end{proof}
\end{lemma}

The key ``computational miracle'' underlying our proof of \Cref{thm:main} is the following. We have been unable to find a suitable generalization beyond the toothbrush case.

\begin{lemma}\label{lem:miracle}
  We have
    \begin{equation}\label{eq:miracle} f^{(2^3, 1^{k-2})}(q) = \qbinom{3+k}{3}_q - q^{k-1} [k+3]_q. \end{equation}

  \begin{proof}
    The the $q$-hook length formula gives
      \[ f^{(2^3, 1^{k-2})}(q) = \frac{[k+4]_q [k+3]_q [k-1]_q}{[3]_q!}. \]
    The right-hand side of Equation~\eqref{eq:miracle} is
      \[ \frac{[k+3]_q [k+2]_q [k+1]_q}{[3]_q!} - \frac{q^{k-1} [k+3]_q [3]_q [2]_q}{[3]_q!}. \]
    The two are equal if and only if
      \[ (1-q^{k+4})(1-q^{k-1}) = (1-q^{k+2})(1-q^{k+1}) - q^{k-1}(1-q^3)(1-q^2), \]
    which may be checked directly.
  \end{proof}
\end{lemma}

We may now restate and prove our main result.

\begin{reptheorem}{thm:main}
  The triple 
  \[
  \left(\IncP^{3+k}(3 \times k), \langle \KPro\rangle, f^{(2^3, 1^{k-2})}(q)\right)
  \]
  exhibits the cyclic sieving phenomenon, and $\KPro^{3+k} = 1$.

  \begin{proof}
    By \Cref{lem:decomp}, we have an equivariant decomposition $\Inc^{3+k}(3 \times k) = \IncP^{3+k}(3 \times k) \sqcup \cE$. The triple $\left(\cE, \langle \KPro\rangle, q^{k-1}[k+3]_q\right)$ exhibits the CSP since $\cE$ is a single orbit of length $k+3$ and
    \begin{align*}
      q^{k-1}[k+3]_q
        &= q^{k-1} + q^k + \cdots + q^{2k+1} \\
        &\equiv 1 + q + \cdots + q^{k+2} \pmod{1-q^{k+3}} \\
        &= [k+3]_q,
    \end{align*}
    using the stabilizer-order criterion for the CSP \cite[p.18]{MR2087303}. By \Cref{cor:Inc_CSP}, this is an instance of \textit{refined cyclic sieving} in the sense of \cite[p.39]{MR3836732}, and it follows that
      \[ \left(\IncP^{3+k}(3 \times k), \langle \KPro\rangle, \qbinom{3+k}{3}_q - q^{k-1} [k+3]_q\right) \]
    exhibits the CSP. The result follows by \Cref{lem:miracle}.
  \end{proof}
\end{reptheorem}

\section{Bijections}\label{sec:bijections}

The original argument for (CSP.1) in \Cref{sec:existing} involves Kaszhdan--Lusztig cellular representations, and no bijective proof is known. By contrast, the argument for (CSP.2) uses a $\maj$-preserving bijection to standard tableaux and direct evaluations at roots of unity. (For more algebraic perspectives on (CSP.2), see \cite{MR3591372,Kim.Rhoades,Patrias.Pechenik.Striker}.) While the argument for (CSP.3) is not bijective, it involves a map to standard tableaux with well-controlled fibers and direct evaluations at roots of unity.

Our proof of \Cref{prop:main} does not produce a single ``natural'' bijection between $\IncP^{3+k}(3 \times k)$ and $\SYT(2^3, 1^{k-2})$. One may, however, produce a bijection as follows. First identify $\SYT(2^3, 1^{k-2})$ with the $3$-element subsets of $\{2, 3, \ldots, k+4\}$ consisting of the entries in the second column. The collection $\cE'$ of ``exceptional'' subsets not of this form consists of
  \[ \{2, 3, 4\}, \ldots, \{2, 3, k+4\}, \{2, 4, 5\}, \{3, 4, 5\}, \]
or $k+3$ in all. Hence we have a bijection $\SYT(2^3, 1^{k-2}) \to \binom{[2, k+4]}{3} - \cE'$. Now \Cref{lem:decomp} provides a bijection $\IncP^{3+k}(3 \times k) \to \binom{[k+3]}{3} - \cE$. Write $\cE_-'$ for $\cE'$ with all entries decreased by $1$, and likewise replace $\binom{[2, k+4]}{3}$ with $\binom{[k+3]}{3}$ by decrementing. Pick any bijection on $\cE \cup \cE_-'$ which sends $\cE$ to $\cE_-'$ and extend this bijection to $\binom{[k+3]}{3}$ as the identity elsewhere.

Unfortunately, this construction appears to have almost no useful properties. It would be very interesting to find a ``natural'' bijection proving \Cref{thm:main}. We note that one may compute that no such bijection exists preserving the major index statistic or intertwining the K-evacuation maps.

\section*{Acknowledgements}

The first ideas for this work appeared at the Dynamical Algebraic Combinatorics workshop, held online through the Banff International Research Station (BIRS) in October 2020. The main results were proven at the subsequent BIRS Dynamical Algebraic Combinatorics workshop hosted by the University of British Columbia Okanagan in November 2021. Both workshops were organized by Striker with James Propp, Tom Roby, and Nathan Williams. We are grateful to all the institutions and people involved with making these workshops a success and providing a conducive working environment.

Gaetz acknowledges support from a National Science Foundation Postdoctoral Research Fellowship (DMS-2103121). Pechenik acknowledges support from a Discovery Grant (RGPIN-2021-02391) and Launch Supplement (DGECR-2021-00010) from the Natural Sciences and Engineering Research Council of Canada. Striker acknowledges support from Simons Foundation/SFARI grant (527204, JS).

\bibliography{refs}{}
\bibliographystyle{amsalpha}

\end{document}